\documentclass[11pt,]{amsart}
\usepackage{latexsym}

\newcommand{\eproof}{\mbox{\ }\hfill $\Box$ \par \vskip 10pt}

\newtheorem{Theorem}{Theorem}[section]
\newtheorem{lemma}[Theorem]{Lemma}
\newtheorem{prop}[Theorem]{Proposition}
\newtheorem{rem}[Theorem]{Remark}

\newtheorem{corol}[Theorem]{Corollary}

\def\R{{\mathbb R}}
\def\C{{\mathbb C}}
\def\N{{\mathbb N}}
\def\cal{\mathcal}
\def\Im{{\rm Im\:}}
\def\Re{{\rm Re\:}}

\usepackage{graphicx}
\graphicspath{%
    {converted_graphics/}% inserted by PCTeX
    {/}% inserted by PCTeX
}
\begin{document}

\title [Interior transmission eigenvalues]{Asymptotics of the number of the interior transmission eigenvalues}
\author[V. Petkov]{Vesselin Petkov}
\author[G. Vodev]{Georgi Vodev}

\address {Universit\'e de Bordeaux, Institut de Math\'ematiques de Bordeaux,  351, Cours de la Lib\'eration, 33405  Talence, France}
\email{petkov@math.u-bordeaux1.fr}
\address {Universit\'e de Nantes, D\'epartement de Math\'ematiques, 2, rue de la Houssini\`ere, 44322 Nantes-Cedex, France}
\email{vodev@math.univ-nantes.fr}
\thanks{The first author was partially supported by the ANR project Nosevol BS01019 01}

\date{}

\maketitle

\noindent
\begin{abstract} We prove Weyl asymptotics $N(r) = c r^d + {\mathcal O}_{\epsilon}(r^{d - \kappa + \epsilon})$, $\forall\, 
0<  \epsilon \ll 1$, for the counting function $N(r) = \sharp\{\lambda_j \in \C \setminus \{0\}:\: |\lambda_j| \leq r^2\}$, $r>1$,  of the interior transmission eigenvalues (ITE), $\lambda_j$. Here $d$ denotes the space dimension and $0<\kappa\le 1$ is such that there are no (ITE) in the region $\{\lambda\in \C:\: |{\rm Im}\,\lambda|\ge C(|
{\rm Re}\,\lambda|+1)^{1-\frac{\kappa}{2}}\}$ for some $C>0$.
\end{abstract}
\vspace{0.4cm}
{\bf Key words:} Interior transmission eigenvalues, Weyl formula with remainder, eigenvalue-free regions\\

{\bf MSC:} Primary  35P20, Secondary 35P15, 35P25 
\vspace{0.2cm}

\setcounter{section}{0}
\section{Introduction and statement of results}

Let $\Omega\subset \R^d$, $d\ge 2$, be a bounded, connected domain with a $C^\infty$ smooth boundary $\Gamma=\partial\Omega$. 
A complex number $\lambda\in \C, \lambda \neq 0,$ will be called an interior transmission eigenvalue (ITE) if the following problem has a non-trivial solution:
$$\left\{
\begin{array}{lll}
\left(\nabla c_1(x)\nabla+\lambda n_1(x)\right)u_1=0 &\mbox{in} &\Omega,\\
\left(\nabla c_2(x)\nabla+\lambda n_2(x)\right)u_2=0 &\mbox{in} &\Omega,\\
u_1=u_2,\,\,\, c_1\partial_\nu u_1=c_2\partial_\nu u_2& \mbox{on}& \Gamma,
\end{array}
\right.
\eqno{(1.1)}
$$
where $\nu$ denotes the exterior Euclidean unit normal to $\Gamma$, $c_j,n_j\in C^\infty(\overline\Omega)$, 
$j=1,2$ are strictly positive real-valued functions. The spectral problem for (ITE) is related to a non self-adjoint operator ${\cal A}$ 
(see Section 3) and in the isotropic case $c_1(x) = c_2(x) = 1$ the boundary problem (1.1) is not parameter-elliptic. 
For these reasons many well-known techniques developed for self-adjoint operators or for parameter-elliptic boundary problems are not applicable. 
The positive (ITE) are related to the inverse scattering problems . More precisely,  if $\lambda = k^2$ is a real (ITE), 
then the far-field operator $F(\lambda): L^2({\mathbb S}^{d-1}) \longrightarrow L^2({\mathbb S}^{n-1})$ 
with kernel the scattering amplitude $s(k, \theta, \omega)$ is not injective and its range is not  dense.  This is crucial for the so-called linear sampling method (see \cite{kn:CK}, \cite{kn:AH}) which works if we avoid the real (ITE).  For this reason the problem of the existence 
and the discreteness of (ITE) draw the attention of many authors (see the survey \cite{kn:CH} for a comprehensive review and a more complete list of references). Secondly, it was proved that we can determine the (ITE) 
from the far-field operator. Finally, it was established that in some cases the knowledge of all complex (ITE) determines the index of refraction 
 of the scattering obstacle (see \cite{kn:CH}, \cite{kn:CL}). This explains the increasing interest toward (ITE) and the fact that 
 a lot of papers concerning the existence and the spectral properties of (ITE) in relation with the inverse scattering problems of reconstruction have been recently published. 
 
 On the other hand, the analysis of the (ITE) leads to some interesting and difficult mathematical spectral problems for non self-adjoint operators. These problems are connected with two {\bf major questions}:\\
 
 (A) Describe the eigenvalue-free regions in the complex plane.\\
 
 (B) Find a Weyl asymptotic of the counting function of the eigenvalues.\\
 
 In contrast to the case of self-adjoint operators these questions are much more difficult and there are no general results. 
 As far as the Weyl asymptotics are concerned, one may study the leading term of the counting function and one can search an optimal remainder.  On the other hand, even in the case of boundary problems for non self-adjoint operators which are parameter-elliptic, the Weyl asymptotics in the literature concern mainly the leading term (see \cite{kn:BK} for some results in this direction for non self-adjoint operators).
 
The question (A) has been investigated by the second author in \cite{kn:V} (see also \cite{ kn:H} for a weaker result) and the result in \cite{kn:V} plays an important role in our analysis. In the present paper our purpose is to study the question (B).  Under some conditions the (ITE) form a discrete set in 
$\C \setminus \{0\}$ and they have as an accumulation point only infinity (see for instance \cite{kn:HK2}, \cite{kn:Sy}).
 Introduce the counting function
$$N(r):=\sharp\left\{\lambda_j \in \C \setminus \{0\}: \lambda_j \:\text{is (ITE)},\, |\lambda_j|\le r^2\right\},\,\, r>1, $$
 where the eigenvalues are counted with their multiplicity (see Section 3 for the precise definition of the multiplicity). Recently,  
 many works concerning the Weyl asymptotics of $N(r)$ have been published  both in the isotropic $(c_1\equiv c_2\equiv  1)$ and anisotropic cases 
 (see \cite{kn:R1}, \cite{kn:DP}, \cite{kn:R2}, \cite{kn:LV1}, \cite{kn:LV3}, \cite{kn:LV4}, \cite{kn:HP}, \cite{kn:F}). In \cite{kn:HP} 
 the case when $\Omega = \{x \in \R^d : \: |x| \leq 1\}$ and $c_1\equiv c_2\equiv 1, \: n_1\equiv 1, n_2 = const \neq 1$ has been investigated and for $d = 1$ a sharp asymptotics of $N(r)$ with remainder ${\cal O} (1)$ has been established. 
 In all other works only the leading term of $N(r)$ was obtained. We should mention that in \cite{kn:LV1} 
 the anisotropic case has been studied and the asymptotics of $N(r)$ with a remainder is stated. However, 
 the proof has a gap and  only the asymptotics with leading term seems to be correct. 
 The isotropic case is more difficult since the boundary problem is not parameter-elliptic and the tools for elliptic boundary problems 
 cannot be applied.  In the isotropic case
 when $n_1(x) \equiv 1$, $n_2(x)>1$, $\forall x \in \bar{\Omega}$, it has been recently 
   established in \cite{kn:F}, \cite{kn:R2}  the asymptotics 
$$N(r) \sim (\tau_1 + \tau_2) r^d,\: \quad r \to +\infty, \eqno(1.2)$$
where $ \tau_1$ and $ \tau_2$ are defined below. It is important to remark that in  \cite{kn:F}, \cite{kn:R2} the analysis is based on the study of some trace class 
operators leading to an asymptotics
$$\sum_j \frac{1}{|\lambda_j|^p + t}= \alpha t^{-1 + \frac{d}{2p}} + o(t^{-1 + \frac{d}{2p}}),\: \quad t \to  +\infty, \eqno{(1.3)}$$
where $p \in \N$ is sufficiently large. Combining this asymptotics with the Tauberian theorem of Hardy and Littlewood, 
one obtains (1.2)  and the remainder is given by the principal part divided by a logarithmic factor. To obtain a sharper remainder one 
could apply  a finer Tauberain theorem (see \cite{kn:M}), but for this purpose it is necessary to establish asymptotics like (1.3) 
with sharper remainder for $t$ lying on certain parabola in $\C$. This, however, seems to be a very difficult problem.

In the present work we follow another approach inspired by the paper \cite{kn:CPV}, where asymptotics have been established 
for the number of the resonances associated to an exterior transmission boundary problem. The purpose is to study the asymptotic 
behavior of $N(r)$ under the condition
$$c_1(x)n_1(x)\neq c_2(x)n_2(x),\quad\forall x\in\Gamma.\eqno{(1.4)}$$
 Our main result is the following

\begin{Theorem}  Assume $(1.4)$ fulfilled. Assume also either the condition
$$c_1(x)=c_2(x),\quad \partial_\nu c_1(x)=\partial_\nu c_2(x),\quad \forall x\in\Gamma,\quad \eqno{(1.5)}$$
or the condition
$$c_1(x)\neq c_2(x),\quad\forall x\in\Gamma.\eqno{(1.6)}$$
Then, the $(ITE)$ form a discrete set in $\C $ and we have the following asymptotics
 $$N(r)=(\tau_1+\tau_2)r^d+{\cal O}_\varepsilon(r^{d-\kappa+\varepsilon}),\:\quad r \to + \infty,\eqno{(1.7)}$$
 for every $0<\varepsilon\ll 1$, where
 $$\tau_j=\frac{\omega_d}{(2\pi)^d}\int_\Omega\left(\frac{n_j(x)}{c_j(x)}\right)^{d/2}dx,$$
$\omega_d$ being the volume of the unit ball in $\R^d$, and $\kappa=\frac{1}{2}$ if $(1.5)$ holds, $\kappa=\frac{2}{5}$ if $(1.6)$ holds. 
 Moreover, if in addition to $(1.6)$ we assume either the condition
$$\frac{n_1(x)}{c_1(x)}\neq \frac{n_2(x)}{c_2(x)},\quad\forall x\in\Gamma,\eqno{(1.8)}$$
or the condition
$$\frac{n_1(x)}{c_1(x)}=\frac{n_2(x)}{c_2(x)},\quad\forall x\in\Gamma,\eqno{(1.9)}$$
then $(1.7)$ holds with $\kappa=\frac{1}{2}$. 
\end{Theorem}

To prove this theorem we use in an essential way the eigenvalue-free regions obtained in \cite{kn:V}. In fact, we prove 
in the present paper a {\bf more general result}  saying that if there are no interior transmission eigenvalues in a region of the form
$$\left\{\lambda\in \C:\: |{\rm Im}\,\lambda|\ge C(|{\rm Re}\,\lambda|+1)^{1-\frac{\kappa}{2}}\right\},\quad C>0,\,\,0<\kappa\le 1,\eqno{(1.10)}$$
then the asymptotics (1.7) with remainder ${\mathcal O}_{\epsilon} (r^{d-\kappa + \epsilon})$ is true. On the other hand, it is proved in \cite{kn:V} that under the assumptions of Theorem 1.1,
we have indeed an eigenvalue-free region (1.10) with $\kappa$ replaced by $\kappa-\epsilon$, where  $\kappa$ is given 
by Theorem 1.1. Note that the parametrix construction of the Dirichlet-to-Neumann map in \cite{kn:Sj}, Section 11, suggests that 
for strictly concave domains there are reasons to believe that (1.10) is true with $\kappa=\frac{2}{3}$. 
Hence in this case our argument will imply immediately a better bound of the remainder in (1.7). 
It is also worth noticing that if we have an eigenvalue-free region of the form (1.10) with $\kappa=1$,
we get asymptotics with an almost optimal remainder term ${\cal O}_\varepsilon(r^{d-1+\varepsilon})$ in (1.7).
However, the existence of such an eigenvalue-free region is a very difficult open problem. Nevertheless, according to our result, the problem of bounding the remainder in the Weyl formula for the (ITE) is reduced to that of getting an eigenvalue-free region in $\C$, and a larger eigenvalue-free region yields a sharper bound for the remainder. To our best knowledge, it seems that our paper is the first one where such a relationship is established.

For reader's convenience, in what follows in this section we will discuss the main steps in the proof of Theorem 1.1. The starting point of our argument is a trace formula (see Section 3 and (3.5)) which allows us to relate the number of the (ITE) with the number 
of the eigenvalues, $\nu_j$, of two self-adjoint operators for which the Weyl asymptotics are known to hold, together with a trace of an operator given by an integral involving a meromorphic operator-valued function, $T(\lambda)$, and its inverse $T^{-1}(\lambda)$ (see formula (3.6)). The main problem to deal with is to estimate the trace of this integral and it yields the bound 
${\mathcal O}_{\epsilon}(r^{d - \kappa + \epsilon})$ of the remainder. 
We apply this formula to obtain an asymptotics for the difference $N(r)-N(r/\sqrt{2}),\: r \to + \infty$, and by a standard 
argument it is easy to see that this is sufficient to prove (1.7). 

Since it is more convenient to
work in the semi-classical setting, we reduce our problem to a semi-classical one by introducing a small parameter 
$h=\frac{\sqrt{2}}{r},\: r \gg 1$. Thus we are going to count
the number of  points $\{z_k\}$, $\frac{z_k}{h^2}$ being an (ITE), in a region of the form 
$$\left\{z\in\C:\,1-A h^{\kappa-\epsilon} \le |{\rm Re}\,z|\le 2+A h^{\kappa-\epsilon} ,\,|{\rm Im}\,z|\le h^{\kappa-\epsilon}
\right\},\quad A > 0,$$
provided
we have an eigenvalue-free region (1.10) with $\kappa-\epsilon$ in place of $\kappa$ (see Proposition 3.7). This requires to make a change of variables $\lambda=z/h^2$ in the trace formula (3.6) and to study the behavior 
of the integral term  when $ 0 < h \leq h_0(\epsilon)$ and 
$$ z \in Z= \{z  \in \C:1/2 <|\Re z|<3,|\Im z|<1\}.$$
Next we construct a meromorphic function $g_h(z)$ with poles among the points $\{h^2\nu_j\}$
and such that if an (ITE), $\lambda_k$, does not belong to the set $\{\nu_j\}$, then $h^2\lambda_k$ is a zero of $g_h(z)$
and the multiplicities of the corresponding zeros of $g_h(z)$ and (ITE) agree.  It should be mentioned that the construction of 
the function $g_h(z)$ is not trivial and it requires to build a semi-classical parametrix for 
the corresponding Dirichlet-to-Neumann map  ${\mathcal N}(z, h)$ in the elliptic zone.
This is carried out in Section 2 by using the parametrix construction in \cite{kn:V}.

The estimate of the remainder is reduced to that of the integral
$$\frac{1}{2 \pi i}\int _{\gamma_0}\frac{d}{dz}\det \log g_h(z) dz,\eqno{(1.11)}$$
where $\gamma_0 \subset Z$ is a suitable closed contour chosen so that on $
\gamma_0$ we have neither zeros nor poles of $g_h(z)$. 
The main property of the function $g_h(z)$ is the estimate 
$$\log|g_h(z)|\le C_\epsilon h^{1-d-\epsilon},\quad\forall\,
0<\epsilon\ll 1,$$ 
provided the distance between $z$ and the set $\{h^2\nu_j\}$ is greater than $h^M$, $ M>0$ being arbitrary (see Lemma 3.4).
This estimate plays a crucial role in the estimate of (1.11). Next in Lemma 3.5 we show that
for $z \in Z,\: |\Im z |\geq h^{\kappa-\epsilon},$ we also have
 $$\log\frac{1}{|g_h(z)|}\le C_\varepsilon h^{1-d-\varepsilon},\quad\forall\,0<\varepsilon\ll 1.$$
 Moreover, the function $\log g_h(z)$ is holomorphic in $z\in Z$, $|{\rm Im}\,z|\ge h^{\kappa-\epsilon}$ and satisfies the bound
 $$\left|\frac{d}{dz}\log g_h(z)\right|  \le \frac{C_\epsilon h^{1-d-2\epsilon}}{|{\rm Im} z|}\eqno{(1.12)}$$
 in the domain
 $$W:=\left\{z\in \C:\frac{2}{3}\le|{\rm Re}\,z|\le \frac{5}{2},\,\, 2h^{\kappa-\epsilon}\le |{\rm Im}\,z|\le \frac{1}{2}\right\}.$$
The next step consists of choosing a closed contour $\gamma_0 = \gamma_1 \cup \gamma_3 \cup \gamma_2 \cup \gamma_4$, where $\gamma_3 \subset W,\:\gamma_4 \subset W$ are linear segments parallel to the real axis.  For the integrals over $\gamma_j, \: j = 3, 4,$ we apply (1.12) and one gets
 $$\left|\int_{\gamma_j}
\frac{d}{dz}\log g_h(z)dz\right|\le
C_\epsilon h^{1-d-3\epsilon},\quad j = 3, 4.\eqno{(1.13)}$$
We take $\gamma_j = [w_j^{-}, \tilde{w}_j^{-}] \cup \tilde{\gamma}_j \cup [\tilde{w}_j^{+}, w_j^{+}],\: j = 1,2$ with suitable contours $\tilde{\gamma}_j$ (see Section 3 for the notation). The estimates of the imaginary parts of the integrals over $\tilde{\gamma}_j, \: j = 1,2,$ are more delicate since these contours cross the positive real axis and we must avoid the points $\{ h^2 \nu_k\}$. Our argument is similar to the choice of the contour in \cite{kn:CPV} and the details are given in Section 3. The main point is Lemma 3.8, where the contours $\tilde{\gamma}_j$ are constructed so that
 $$\left|{\rm Im}\,\int_{\widetilde\gamma_j}
\frac{d}{dz}\log g_h(z)dz\right|\le C_\epsilon h^{-d+\kappa-2\epsilon}, \:\quad j = 1,2.\eqno{(1.14)}$$
Combining this with (1.13), we obtain the statement of Proposition 3.7 and by scaling we get the asymptotics of $N(r) - N(r/\sqrt{2}).$\\

\noindent
{\bf Acknowledgment.} Thanks are due to the referee for the comments and remarks concerning the initial version of the paper.

\section{Parametrix of the Dirichlet-to-Neumann map in the elliptic zone}

Let $f\in H^1(\Gamma)$ and consider the problem
$$\left\{
\begin{array}{lll}
 \left(P(h)-z\right)u=0&\mbox{in}&\Omega,\\
 u=f&\mbox{on}&\Gamma,
\end{array}
\right.
\eqno{(2.1)}
$$
where 
$$P(h)=-\frac{h^2}{n(x)}\nabla c(x)\nabla,$$
 $0< h \ll 1$, $z\in Z=\{z\in\C:\frac{1}{2}<|{\rm Re}\,z|<3, |{\rm Im}\,z|<1\}$, $c,n\in C^\infty(\overline\Omega)$ 
 being strictly positive functions.
 The Dirichlet-to-Neumann map is defined by
 $${\cal N}(z,h)f:=\gamma{\cal D}_\nu u:H^{m+1}(\Gamma)\to H^m(\Gamma),$$
 where $m\ge 0$, ${\cal D}_\nu=-ih\partial_\nu$ and $\gamma$ denotes the restriction on $\Gamma$. 
 Denote by $G_D$ the Dirichlet self-adjoint realization of the operator
$-n^{-1}\nabla c\nabla$ on the Hilbert space $H=L^2(\Omega,n(x)dx)$. It is well-known that the spectrum of $G_D$
consists of a discrete set of positive eigenvalues which are also poles of the resolvent $(\lambda - G_D)^{-1}$.
Moreover, if $\nu_k\in{\rm spec}\,G_D$, we have 
$$(\lambda-G_D)^{-1}=\frac{\Pi_k}{\lambda-\nu_k}$$
modulo an operator-valued function holomorphic at $\nu_k$, where $\Pi_k$ is a finite rank projection. The multiplicity of
$\nu_k$ is defined as being the rank of $\Pi_k$. Let ${\cal V}(h):=\{\nu_k\in{\rm spec}\,G_D: h^2\nu_k\in Z\}$. 
The following properties of the Dirichlet-to-Neumann map are more or less 
well-known but  we will give a proof for the sake of completeness. 

\begin{lemma} The Dirichlet-to-Neumann map ${\cal N}(z,h)$ is a meromorphic operator-valued function in $z\in Z$ with poles at $h^2\nu_k$, 
$\nu_k\in {\cal V}(h)$. Moreover,
$${\cal N}(z,h)=\frac{\widetilde\Pi_k(h)}{z-h^2\nu_k}\eqno{(2.2)}$$
modulo an operator-valued function holomorphic at $h^2\nu_k$, where $\widetilde\Pi_k(h)$ is of rank $\le{\rm mult}(\nu_k)$.
If $\delta(z,h):=\min\{1,{\rm dist}\{z,{\rm spec}\,h^2G_D\}\}>0$, then we have the bound
$$\left\|{\cal N}(z,h)\right\|_{H^{m+1}(\Gamma)\to H^m(\Gamma)}\le\frac{Ch}{\delta(z,h)},\eqno{(2.3)}$$
where $C>0$ is a constant which may depend on $m$.
\end{lemma}
 
 {\it Proof.} Clearly, there exists an extension operator $E_m:H^{m+1}(\Gamma)\to H^{m+3/2}(\Omega)$ such that $\gamma E_mf=f$ and 
 $E_mf$ is supported near $\Gamma$. If $f\in H^{m+1}(\Gamma)$ and $z/h^2$ does not belong to ${\rm spec}\,G_D$, 
 it is easy to see that the solution $u$ of (2.1) can be expressed by the formula
 $$u=E_mf-(h^2G_D-z)^{-1}(P(h)-z)E_mf.$$
 Hence
  $${\cal N}(z,h)f=\gamma{\cal D}_\nu E_mf-\gamma{\cal D}_\nu (h^2G_D-z)^{-1}(P(h)-z)E_mf.\eqno{(2.4)}$$
 It follows from (2.4) that ${\cal N}(z,h)$ is a meromorphic operator-valued function in $z\in Z$ with poles among the poles of
 $(h^2G_D-z)^{-1}$ and that (2.2) holds with
 $$\widetilde \Pi_k(h)=\gamma{\cal D}_\nu \Pi_k(P(h)-h^2\nu_k)E_m.$$
 This implies ${\rm rank}\,\widetilde\Pi_k(h)\le {\rm rank}\,\Pi_k$ as desired. By (2.4) we also have
 $$\left\|{\cal N}(z,h)f\right\|_{H^m(\Gamma)}\le Ch\|f\|_{H^{m+1}(\Gamma)}$$ 
 $$+Ch\left\|(h^2G_D-z)^{-1}\right\|_{H^{m+3/2}(\Omega)\to H^{m+3/2}(\Omega)}\left\|E_mf\right\|_{H^{m+3/2}(\Omega)}$$
 $$+Ch\left\|(h^2G_D-z)^{-1}\right\|_{H^{m-1/2}(\Omega)\to H^{m+3/2}(\Omega)}\left\|P(h)E_mf\right\|_{H^{m-1/2}(\Omega)}.$$
 Clearly, we have
 $$\left\|E_mf\right\|_{H^{m+3/2}(\Omega)}\le C\|f\|_{H^{m+1}(\Gamma)},$$
 $$\left\|P(h)E_mf\right\|_{H^{m-1/2}(\Omega)}\le Ch^2\left\|E_mf\right\|_{H^{m+3/2}(\Omega)}\le Ch^2\|f\|_{H^{m+1}(\Gamma)}.$$
 On the other hand, the coercive estimate
$$\left\|v\right\|_{H^{s+2}(\Omega)}\le C\left\|G_Dv\right\|_{H^{s}(\Omega)}+C\left\|v\right\|_{H^{s}(\Omega)},\quad\forall
v\in D(G_D)\cap H^{s}(\Omega)$$
implies the bounds
$$\left\|(h^2G_D-z)^{-1}\right\|_{H^{m+3/2}(\Omega)\to H^{m+3/2}(\Omega)}\le 
C\left\|(h^2G_D-z)^{-1}\right\|_{L^{2}(\Omega)\to L^{2}(\Omega)}\le \frac{C}{\delta(z,h)},$$
$$\left\|(h^2G_D-z)^{-1}\right\|_{H^{m-1/2}(\Omega)\to H^{m+3/2}(\Omega)}$$ $$\le 
C\left\|(G_D-i)(h^2G_D-z)^{-1}\right\|_{H^{m+3/2}(\Omega)\to H^{m+3/2}(\Omega)}\left\|(G_D-i)^{-1}\right\|_{H^{m-1/2}(\Omega)
\to H^{m+3/2}(\Omega)}$$
$$\le \frac{\widetilde C}{h^2\delta(z,h)}.$$
Therefore, (2.3) follows from the above estimates and the proof is complete.
\eproof

Let $(x',\xi')$ be coordinates on $T^*\Gamma$ and denote by $r_0(x',\xi')$ the principal symbol of the Laplace-Beltrami operator,
$-\Delta_\Gamma$, on $\Gamma$ equipped with the Riemannian metric induced by the Euclidean one in $\R^d$. It is well-known that
$r_0$ is a polynomial function in $\xi'$, homogeneous of order 2, and $C_2|\xi'|^2\ge r_0(x',\xi')\ge C_1|\xi'|^2$ with constants $C_2>C_1>0$.
Set $m(x)=\frac{n(x)}{c(x)}$. Let $\phi\in C^\infty(\R)$,
$\phi(\sigma)=1$ for $|\sigma|\le 1$, $\phi(\sigma)=0$ for $|\sigma|\ge 2$, and set
$$\chi(x',\xi')=\phi\left(\delta_0r_0(x',\xi')\right),$$
where $0<\delta_0\ll 1$. For $(x',\xi')\in{\rm supp}\,(1-\chi)$, 
introduce the function
$$\rho(x',\xi',z)=i\sqrt{r_0(x',\xi')-\gamma m(x')z}=i\sqrt{r_0}\left(1-z\frac{\gamma m}{r_0}\right)^{1/2}.$$
Since 
$$|z|\frac{\gamma m}{r_0}\le\frac{1}{2},\quad\forall z\in Z,\,(x',\xi')\in{\rm supp}\,(1-\chi),$$
the functions $\rho$ and $\rho^{-1}$ are holomorphic in $z\in Z$ and 
$${\rm Im}\,\rho(x',\xi',z)\ge C\sqrt{r_0(x',\xi')}$$
with some constant $C>0$. 
  In what follows in this section  we will construct a parametrix for the operator ${\cal N}(z,h){\rm Op}_h(1-\chi)$, where
  ${\rm Op}_h(1-\chi)$ denotes the $h-\Psi$DO with symbol $1-\chi$. In fact, this construction is carried out in \cite{kn:V} and here we will
  only recall the main points. First, notice that it suffices to make the construction
 locally and then to glue up all pieces by using a partition of the unity on $\Gamma$.  
Given an arbitrary point $x^0\in\Gamma$, there exists a small neighborhood ${\cal O}(x^0)\subset\overline\Omega$ of $x^0$
and local normal coordinates $(x_1,x')\in {\cal O}(x^0)$ such that $x^0=(0,0)$, $\Gamma\cap {\cal O}(x^0)$ is defined by $x_1=0$,
$x'$ being coordinates in $\Gamma\cap {\cal O}(x^0)$, $x_1>0$ in $\Omega\cap {\cal O}(x^0)$, and in these coordinates the operator
$${\cal P}(z,h)=-\frac{h^2}{c(x)}\nabla c(x)\nabla-z\frac{n(x)}{c(x)}$$
can be written in the form
$${\cal P}(z,h)={\cal D}_{x_1}^2+r(x,{\cal D}_{x'})-zm(x) + hq(x,{\cal D}_x)+h^2\widetilde q(x).$$
Here we have set ${\cal D}_{x_1}=-ih\partial_{x_1}$, ${\cal D}_{x'}=-ih\partial_{x'}$, $r(x,\xi')=\langle R(x)\xi',\xi'\rangle$, 
$R=(R_{ij})$ being a symmetric $(d-1)\times(d-1)$ matrix-valued function with smooth real-valued entries, $q(x,\xi)=\langle q(x),\xi\rangle$, 
$q(x)$ and $\widetilde q(x)$ being smooth functions. 
Moreover, we have $r(0,x',\xi')=r_0(x',\xi')$, $r_0(x', \xi')$ being the principal symbol of $-\Delta_\Gamma$ written in 
the coordinates $(x',\xi')$.
Let $\psi(x')\in C_0^\infty(\Gamma\cap {\cal O}(x^0))$, $\psi=1$ in a neighborhood of $x^0$. In \cite{kn:V}, 
it was constructed a parametrix, $\widetilde u_\psi$, 
of (2.1) satisfying the condition $\widetilde u_\psi|_{x_1=0}={\rm Op}_h(1-\chi)\psi f$ and having the form
$$\widetilde u_\psi(x)=(2\pi h)^{-d+1}\int\int e^{\frac{i}{h}\varphi(x,y',\xi',z)}\phi\left(\frac{x_1}{\delta_1}
\right)a(x,\xi',z,h)f(y')dy'd\xi',$$
where $\phi$ is as above and $\delta_1>0$ is a small constant independent of $x,\xi',h,z$. The phase $\varphi$ is a
complex-valued function such that 
$$\varphi|_{x_1=0}=-\langle x'-y',\xi'\rangle, \:\partial_{x_1}\varphi|_{x_1=0}=\rho,\:{\rm Im}\,\varphi\ge x_1{\rm Im}\,\rho/2,$$
 and the amplitude $a$ satisfies $a|_{x_1=0}=\psi(x')(1-\chi(x',\xi'))$. More generally, the functions $\varphi$ and $a$ are of the form
$$\varphi=-\langle x'-y',\xi'\rangle+\sum_{k=1}^{N-1} x_1^k\varphi_k(x',\xi',z)=-\langle x'-y',\xi'\rangle+\widetilde\varphi,$$
$$a=\sum_{k=0}^{N-1} \sum_{j=0}^{N-1} x_1^kh^ja_{k,j}(x',\xi',z),$$
$N\gg 1$ being an arbitrary integer. The phase $\varphi$ satisfies the eikonal equation mod ${\cal O}(x_1^N)$:
$$\left(\partial_{x_1}\varphi\right)^2+r(x,\nabla_{x'}\varphi)-m(x)z=x_1^N\Psi_N(x,\xi',z)\eqno{(2.5)}$$
 and $a$ satisfies the equation
$$e^{-\frac{i}{h}\varphi}{\cal P}(z,h)e^{\frac{i}{h}\varphi}a=x_1^NA_N(x,\xi',z,h)+h^NB_N(x,\xi',z,h),\eqno{(2.6)}$$
where $\Psi_N$, $A_N$ and $B_N$ are smooth functions. It was shown in Section 4 of \cite{kn:V} that  $a_{k,j}\in S^{-j}$, $j\ge 0$, $k\ge 1$, 
$\partial_{x_1}^kA_N\in S^2$, $\partial_{x_1}^kB_N\in S^{1-N}$, $k\ge 0$, uniformly in $z\in Z$ and $0<x_1\le\delta_1$. Recall that
$S^k$ are the spaces of all functions $a\in C^\infty(T^*\Gamma�)$ satisfying the estimates
$$\left|\partial_{x'}^\alpha\partial_{\xi'}^\beta a(x',\xi')\right|\le C_{\alpha,\beta}\langle\xi'\rangle^{k-|\beta|},\: \langle \xi'\rangle = (1 + |\xi'|^2)^{1/2}$$
for all multi-indices $\alpha$ and $\beta$.
Moreover, the functions $a_{k,j}$, $A_N$, $B_N$ are polynomials in $\rho$, $\rho^{-1}$ and $z$, and therefore they are holomorphic in $z\in Z$.
As in \cite{kn:V}, it is easy to see that
$${\cal P}(z,h)\widetilde u_\psi={\rm Op}_h(p_\psi)f,$$
where the function
$$p_\psi=e^{\frac{i}{h}\langle x',\xi'\rangle}\left[{\cal P}(z,h),\phi\left(\frac{x_1}{\delta_1}\right)\right]e^{-\frac{i}{h}
\langle x',\xi'\rangle}e^{\frac{i}{h}\widetilde\varphi}a
+e^{\frac{i}{h}\widetilde\varphi}\phi\left(\frac{x_1}{\delta_1}\right)\left(x_1^NA_N+h^NB_N\right)$$
is holomorphic in $z$ and satisfies the bounds
$$|\partial_x^\alpha p_\psi|\le C_{\alpha,N}\left(\frac{h}{\langle\xi'\rangle}\right)^{N-\ell-|\alpha|}\quad
\mbox{for}\quad|\alpha|\le N-\ell\eqno{(2.7)}$$
with some $\ell$ independent of $N$ and $\alpha$.
The parametrix, $\widetilde{\cal N}_\psi(z,h)$, of the operator ${\cal N}(z,h){\rm Op}_h(1-\chi)\psi$ is defined by
$${\cal D}_{x_1}\widetilde u_\psi|_{x_1=0}=\widetilde{\cal N}_\psi(z,h)f={\rm Op}_h(\eta_\psi)f,$$
where 
$$\eta_\psi=a\frac{\partial\varphi}{\partial x_1}|_{x_1=0}-ih\frac{\partial a}{\partial x_1}|_{x_1=0}=\psi(1-\chi)\rho-ih\sum_{j=0}^{N-1}
h^ja_{1,j},$$
$$a_{1,0}=-\frac{i}{2} q(0,x',1,\xi'/\rho)\psi-\frac{1}{2\rho}\langle R(0,x')\xi',\nabla_{x'}\psi(x')\rangle.$$
Since
$$\frac{1}{\rho}=\frac{1}{i\sqrt{r_0}}\left(1-z\frac{\gamma m}{r_0}\right)^{-1/2}
=\frac{1}{i\sqrt{r_0}}+{\cal O}\left(\langle\xi'\rangle^{-3}\right),$$
 we deduce that mod $S^{-2}$ the function $a_{1,0}$ is given by the expression
$$a_{1,0}=-\frac{1}{2} q(0,x',1,\xi'/\sqrt{r_0})\psi+\frac{i}{2\sqrt{r_0}}\langle R(0,x')\xi',\nabla_{x'}\psi(x')\rangle$$
 $$=\left\langle \frac{i\nabla_{x'}c(0,x')}{2c(0,x')},\frac{\xi'}{\sqrt{r_0}}\right\rangle\psi+
 \frac{i\partial_{x_1}c(0,x')}{2c(0,x')}\psi+q_0(x',\xi')\eqno{(2.8)}$$
 with some function $q_0\in S^0$ independent of the functions $c$ and $n$. 
 
Let $\{\psi_j\}_{j=1}^J$ be a partition of the unity on $\Gamma$. Set
$$p=\sum_{j=1}^Jp_{\psi_j},\quad \eta=\sum_{j=1}^J\eta_{\psi_j},\quad \widetilde u=\sum_{j=1}^J\widetilde u_{\psi_j}.$$
The operator
 $$\widetilde{\cal N}(z,h)=\sum_{j=1}^J \widetilde{\cal N}_{\psi_j}(z,h)={\rm Op}_h(\eta)$$
  is an $h-\Psi$DO on $\Gamma$ with a principal symbol $\rho(1-\chi)$, holomorphic in $z\in Z$. 
Let $u_{\psi_j}$ be the solution of (2.1) with $u_\psi|_{\Gamma}={\rm Op}_h(1-\chi)\psi f$. 
Then $u=\sum_{j=1}^Ju_{\psi_j}$ is the solution of (2.1) with $u|_{\Gamma}={\rm Op}_h(1-\chi)f$.  
Moreover, it is easy to see that, if $z/h^2$ does not belong to spec$\,G_D$, we have 
$$u=\widetilde u-\left(h^2G_D-z\right)^{-1}\frac{c}{n}{\cal P}(z,h)\widetilde u$$
which yields the identity
$${\cal N}(z,h){\rm Op}_h(1-\chi)f=\widetilde{\cal N}(z,h)f-\gamma{\cal D}_\nu\left(h^2G_D-z\right)^{-1}\frac{c}{n}{\rm Op}_h(p)f.\eqno{(2.9)}$$
It follows from (2.7) that if $N$ is taken large enough, the operator
$$F(z,h):={\cal N}(z,h)-\widetilde{\cal N}(z,h)={\cal N}(z,h){\rm Op}_h(\chi)
-\gamma{\cal D}_\nu\left(h^2G_D-z\right)^{-1}\frac{c}{n}{\rm Op}_h(p)$$
is meromorphic with values in the space of trace class operators on $L^2(\Gamma)$. Let $\mu_j(F)$ be the characteristic values of $F$.
Recall that $\mu_j(F)$ are defined as being the eigenvalues of the self-adjoint operator $(F^*F)^{1/2}$. 

\begin{lemma} If $z/h^2$ does not belong to ${\rm spec}\,G_D$, then for every integer $0\le m\le N/4$ we have the bound
$$\mu_j(F(z,h))\le \frac{C}{\delta(z,h)}\left(hj^{1/(d-1)}\right)^{-2m},\: \forall j, \eqno{(2.10)}$$
where the constant $C>0$ depends on $m$ and $N$ but is independent of $z$, $h$, $j$, and $\delta(z,h)$ is defined in Lemma $2.1.$
\end{lemma}
 
 {\it Proof.} We will use the well-known fact that the characteristic values of the Laplace-Beltrami operator 
 on a compact Riemannian manifold without boundary
 (in our case $\Gamma$, dim\,$\Gamma=d-1$) satisfy
 $$\mu_j\left((1-\Delta_\Gamma)^{-m}\right)\le C_mj^{-2m/(d-1)},\: \forall j, \eqno{(2.11)}$$
 for every integer $m\ge 0$. On the other hand, by using the trace theorem and Lemma 2.1, we obtain
 $$\left\|F(z,h)\right\|_{L^2(\Gamma)\to H^{2m}(\Gamma)}\le \left\|{\cal N}(z,h)\right\|_{H^{2m+1}(\Gamma)\to H^{2m}(\Gamma)}
 \left\|{\rm Op}_h(\chi)\right\|_{L^2(\Gamma)\to H^{2m+1}(\Gamma)}$$
  $$+Ch\left\|\left(h^2G_D-z\right)^{-1}\right\|_{H^{2m+3/2}(\Omega)\to H^{2m+3/2}(\Omega)}
  \left\|{\rm Op}_h(p)\right\|_{L^2(\Gamma)\to H^{2m+3/2}(\Omega)}$$
  $$\le \frac{Ch}{\delta(z,h)}\left\|{\rm Op}_h(\chi)\right\|_{L^2(\Gamma)\to H^{2m+1}(\Gamma)}
  +\frac{Ch}{\delta(z,h)}\left\|{\rm Op}_h(p)\right\|_{L^2(\Gamma)\to H^{2m+3/2}(\Omega)}.$$
  Since the function $\chi$ is compactly supported, we have the bound 
  $$\left\|{\rm Op}_h(\chi)\right\|_{L^2(\Gamma)\to H^{2m+1}(\Gamma)}\le C_mh^{-2m-1}.\eqno{(2.12)}$$
  In view of (2.7) we also have
 $$ \left\|{\rm Op}_h(p)\right\|_{L^2(\Gamma)\to H^{2m+3/2}(\Omega)}\le C_{m,N}h^{N-2m-\ell_1}\eqno{(2.13)}$$
 with some $\ell_1$ independent of $m$ and $N$, provided $0\le m\le N/4$ and  $N$ being large enough. By (2.12) and (2.13) we conclude
  $$\left\|F(z,h)\right\|_{L^2(\Gamma)\to H^{2m}(\Gamma)}\le \frac{C_mh^{-2m}}{\delta(z,h)}.\eqno{(2.14)}$$
  Clearly, (2.10) follows from (2.11) and (2.14) and the proof is complete.
  \eproof

\section{Analysis of the transmission eigenvalues}

For $\lambda\in\C\setminus \{0\}$ define the operator  $R(\lambda)v=u$, where $u=(u_1,u_2)$ and $v=(v_1,v_2)$ solve the problem
$$\left\{
\begin{array}{lll}
\left(-\frac{1}{n_1(x)}\nabla c_1(x)\nabla-\lambda\right)u_1=v_1 &\mbox{in} &\Omega,\\
\left(-\frac{1}{n_2(x)}\nabla c_2(x)\nabla-\lambda\right)u_2=v_2 &\mbox{in} &\Omega,\\
u_1=u_2,\,\,\, c_1\partial_\nu u_1=c_2\partial_\nu u_2& \mbox{on}& \Gamma.
\end{array}
\right.
\eqno{(3.1)}
$$
Denote by $G_D^{(j)}$, $j=1,2$, the Dirichlet self-adjoint realization of the operator
$-n_j^{-1}\nabla c_j\nabla$ on the Hilbert space $H_j=L^2(\Omega,n_j(x)dx)$. Set
${\cal H}=H_1\oplus H_2$ and  define also the operators $K_j(\lambda)f=u$, where $u$ is the solution of the problem
$$\left\{
\begin{array}{lll}
\left(-\frac{1}{n_j(x)}\nabla c_j(x)\nabla-\lambda\right)u=0&\mbox{in} &\Omega,\\
u=f& \mbox{on}& \Gamma.
\end{array}
\right.
\eqno{(3.2)}
$$
Differentiating this equation with respect to $\lambda$, one obtains easily the identity
$$\frac{dK_j(\lambda)}{d\lambda}=(G_D^{(j)}-\lambda)^{-1}K_j(\lambda).\eqno{(3.3)}$$
Introduce the operator 
$$T(\lambda) :=c_1\gamma\partial_\nu K_1(\lambda)-c_2\gamma\partial_\nu K_2(\lambda).$$

\begin{prop} If $T(\lambda)^{-1}$ is a meromorphic operator-valued function with residue of finite rank, the same is true for  
$R(\lambda)$ and we have the formula
$$R(\lambda)=
\left(
\begin{array}{ll}
R_{11}(\lambda),&R_{12}(\lambda) \\
R_{21}(\lambda),&R_{22}(\lambda)
\end{array}
\right):{\cal H}\to {\cal H},
\eqno{(3.4)}$$ 
where
$$R_{11}(\lambda)=(G_D^{(1)}-\lambda)^{-1}-K_1(\lambda)T(\lambda)^{-1}c_1\gamma\partial_\nu(G_D^{(1)}-\lambda)^{-1},$$
$$R_{22}(\lambda)=(G_D^{(2)}-\lambda)^{-1}+K_2(\lambda)T(\lambda)^{-1}c_1\gamma\partial_\nu(G_D^{(2)}-\lambda)^{-1},$$
$$R_{12}(\lambda)=K_1(\lambda)T(\lambda)^{-1}c_1\gamma\partial_\nu(G_D^{(2)}-\lambda)^{-1},$$
$$R_{21}(\lambda)=-K_2(\lambda)T(\lambda)^{-1}c_2\gamma\partial_\nu(G_D^{(1)}-\lambda)^{-1}.$$
Moreover, if $\gamma_0\subset \C$ is a simple closed positively oriented curve which avoids the eigenvalues of 
$G_D^{(j)}$, $j=1,2$, as well as the poles of $T(\lambda)^{-1}$, then we have the identity
$$-{\rm tr}_{{\cal H}}\, (2\pi i)^{-1}\int_{\gamma_0}R(\lambda)d\lambda+\sum_{j=1}^2{\rm tr}_{H_j}\, (2\pi i)^{-1}\int_{\gamma_0}
(G_D^{(j)}-\lambda)^{-1}d\lambda$$ $$={\rm tr}_{L^2(\Gamma)}\, (2\pi i)^{-1}\int_{\gamma_0}T(\lambda)^{-1}\frac{dT(\lambda)}{d\lambda}d\lambda.
\eqno{(3.5)}$$
\end{prop}

{\it Proof.} Clearly, if $(u_j,v_j)$ satisfies (3.1) and $\lambda$ does not belong to ${\rm spec}\,G_D^{(1)}\cup {\rm spec}\,G_D^{(2)}$, we have
$$u_j=(G_D^{(j)}-\lambda)^{-1}v_j+K_j(\lambda)f,$$
where $f=\gamma u_1=\gamma u_2$. The boundary condition in (3.1) implies the identity
$$0 = c_1 \partial_{\nu} u_1 - c_2 \partial_{\nu} u_2 =T(\lambda)f+c_1\gamma\partial_{\nu}(G_D^{(1)}-\lambda)^{-1}v_1- 
c_2\gamma\partial_{\nu}(G_D^{(2)}-\lambda)^{-1}v_2.$$
Hence
$$u_j=(G_D^{(j)}-\lambda)^{-1}v_j-K_j(\lambda)T(\lambda)^{-1}\left(c_1\gamma\partial_{\nu}(G_D^{(1)}-\lambda)^{-1}v_1- 
c_2\gamma\partial_{\nu}(G_D^{(2)}-\lambda)^{-1}v_2\right)$$
which clearly implies (3.4). Moreover, if $T(\lambda)^{-1}$ is meromorphic, so are the operators $R_{ij}(\lambda)$, and by (3.4) the operator
$R(\lambda)$ is meromorphic, too. Using (3.3) and the cyclicity of the trace (see Lemma 2.2 of \cite{kn:SjV}), we get
$${\rm tr}_{{\cal H}}\, (2\pi i)^{-1}\int_{\gamma_0}R(\lambda)d\lambda={\rm tr}_{H_1}\, (2\pi i)^{-1}\int_{\gamma_0}R_{11}(\lambda)d\lambda
+{\rm tr}_{H_2}\, (2\pi i)^{-1}\int_{\gamma_0}R_{22}(\lambda)d\lambda$$
 $$={\rm tr}_{H_1}\, (2\pi i)^{-1}\int_{\gamma_0}(G_D^{(1)}-\lambda)^{-1}d\lambda-{\rm tr}_{H_1}\, (2\pi i)^{-1}\int_{\gamma_0}
 K_1(\lambda)T(\lambda)^{-1}c_1\gamma\partial_\nu(G_D^{(1)}-\lambda)^{-1}d\lambda$$ $$
+{\rm tr}_{H_2}\, (2\pi i)^{-1}\int_{\gamma_0}(G_D^{(2)}-\lambda)^{-1}d\lambda
+{\rm tr}_{H_2}\, (2\pi i)^{-1}\int_{\gamma_0}K_2(\lambda)T(\lambda)^{-1}c_2\gamma\partial_\nu(G_D^{(2)}-\lambda)^{-1}d\lambda$$
$$={\rm tr}_{H_1}\, (2\pi i)^{-1}\int_{\gamma_0}(G_D^{(1)}-\lambda)^{-1}d\lambda-{\rm tr}_{L^2(\Gamma)}\, (2\pi i)^{-1}\int_{\gamma_0}
 T(\lambda)^{-1}c_1\gamma\partial_\nu(G_D^{(1)}-\lambda)^{-1}K_1(\lambda)d\lambda$$ $$
+{\rm tr}_{H_2}\, (2\pi i)^{-1}\int_{\gamma_0}(G_D^{(2)}-\lambda)^{-1}d\lambda
+{\rm tr}_{L^2(\Gamma)}\, (2\pi i)^{-1}\int_{\gamma_0}T(\lambda)^{-1}c_2\gamma\partial_\nu(G_D^{(2)}-\lambda)^{-1}K_2(\lambda)d\lambda$$ 
$$={\rm tr}_{H_1}\, (2\pi i)^{-1}\int_{\gamma_0}(G_D^{(1)}-\lambda)^{-1}d\lambda+{\rm tr}_{H_2}\, 
(2\pi i)^{-1}\int_{\gamma_0}(G_D^{(2)}-\lambda)^{-1}d\lambda$$
$$-{\rm tr}_{L^2(\Gamma)}\, (2\pi i)^{-1}\int_{\gamma_0}
 T(\lambda)^{-1}c_1\frac{d\gamma\partial_\nu K_1(\lambda)}{d\lambda}d\lambda 
 +{\rm tr}_{L^2(\Gamma)}\, (2\pi i)^{-1}\int_{\gamma_0}T(\lambda)^{-1}c_2\frac{d\gamma\partial_\nu K_2(\lambda)}{d\lambda}d\lambda$$ 
which implies (3.5).
\eproof

If $R(\lambda)$ is a meromorphic operator-valued function with residue of finite rank, we define the multiplicity of a pole
$\lambda_k\in\C$ of $R(\lambda)$ by
$${\rm mult}\,(\lambda_k)={\rm rank}\,(2\pi i)^{-1}\int_{|\lambda-\lambda_k|=\varepsilon}R(\lambda)d\lambda,\quad 0<\varepsilon\ll 1.$$
Let the curve $\gamma_0$ be as in Proposition 3.1 and denote by $M_{\gamma_0}$ and $M_{\gamma_0}^{(j)}$, $j=1,2$, the number (counted with the
multiplicity) of the poles
of $R(\lambda)$ and the eigenvalues of $G_D^{(j)}$, respectively, in the interior of $\gamma_0$. Proposition 3.1 implies the
following

\begin{corol} We have the identity
$$M_{\gamma_0}=M_{\gamma_0}^{(1)}+M_{\gamma_0}^{(2)}+{\rm tr}_{L^2(\Gamma)}\, 
(2\pi i)^{-1}\int_{\gamma_0}T(\lambda)^{-1}\frac{dT(\lambda)}{d\lambda}d\lambda.
\eqno{(3.6)}$$
\end{corol}

{\it Proof.} It is easy to see that $R(\lambda)=({\cal A}-\lambda)^{-1}$, where the operator ${\cal A}$ is defined by
$${\cal A}
\begin{pmatrix}
u_1\\
u_2
\end{pmatrix}
=\left(
\begin{array}{l}
-\frac{1}{n_1(x)}\nabla c_1(x)\nabla u_1\\
-\frac{1}{n_2(x)}\nabla c_2(x)\nabla u_2
\end{array}
\right)
$$
with domain 
$$D({\cal A})=\left\{(u_1,u_2)\in {\cal H}: \nabla c_1(x)\nabla u_1\in L^2(\Omega),\,\nabla c_2(x)\nabla u_2\in L^2(\Omega),\right.$$ 
$$\left.\gamma u_1=\gamma u_2,\, c_1\gamma\partial_\nu u_1=c_2\gamma\partial_\nu u_2\right\}.$$
Hence the finite-rank operator
$$-(2\pi i)^{-1}\int_{|\lambda-\lambda_k|=\varepsilon}R(\lambda)d\lambda=(2\pi i)^{-1}\int_{|\lambda-\lambda_k|=\varepsilon}(
{\lambda-\cal A})^{-1}d\lambda$$
is in fact a projection (e.g. see \cite{kn:K}), and therefore the rank coincides with the trace. Thus, (3.6) follows from (3.5).
\eproof

Let $z$ and $h$ be as in the previous section and denote by ${\cal N}_j$, $\widetilde{\cal N}_j$, $F_j$, $j=1,2$, 
the operators defined by replacing in the definition of ${\cal N}$, $\widetilde{\cal N}$, $F$ introduced in Section 2 the pair
$(c,n)$ by $(c_j,n_j)$. Clearly, we have the relationship
$$hT(z/h^2)=c_1{\cal N}_1(z,h)-c_2{\cal N}_2(z,h)$$ $$=c_1F_1(z,h)-c_2F_2(z,h)+c_1\widetilde{\cal N}_1(z,h)-
c_2\widetilde{\cal N}_2(z,h).\eqno{(3.7)}$$
In what follows $H^s_h$ will denote the Sobolev space $H^s(\Gamma)$ equipped with the semi-classical norm.

\begin{lemma} There exist an invertible, bounded operator $E(z,h):H_h^s\to H_h^{s+k}=O(1)$, with an inverse 
$E(z,h)^{-1}:H_h^s\to H_h^{s-k}=O(1)$, $\forall s\in \R$, and trace class operators $L_l(z,h)$ and $L_r(z,h)$ such that
$$E(z,h)\left(c_1\widetilde{\cal N}_1(z,h)-c_2\widetilde{\cal N}_2(z,h)\right)=I +L_l(z,h),\eqno{(3.8)}$$
$$\left(c_1\widetilde{\cal N}_1(z,h)-c_2\widetilde{\cal N}_2(z,h)\right)E(z,h)=I +L_r(z,h),\eqno{(3.9)}$$
where $k=-1$ if $(1.5)$ holds, $k=1$ if $(1.6)$ holds. Moreover, the operators $E, E^{-1}, L_l, L_r$ 
are holomorphic with respect to $z$ for $ z \in Z.$
\end{lemma}

{\it Proof.} Set $m_j=\frac{n_j}{c_j}$, $\rho_j=i\sqrt{r_0-z\gamma m_j},\: j = 1,2,$ and let the real-valued function 
$\chi,\: 0 \leq \chi \leq 1$ be as in Section 2, with a sufficiently large support. It follows from the parametrix construction in 
Section 2 that $c_1\widetilde{\cal N}_1-c_2\widetilde{\cal N}_2={\rm Op}_h(b)$
with a symbol $b=\sum_{j=0}^Nh^j\,b_j$, where $b_j\in S^{1-j}$ are holomorphic in $z\in Z$, and
$$b_0=(c_1\rho_1-c_2\rho_2)(1-\chi).$$
Let $\chi_0\in C^\infty(T^*\Gamma)$ be a real-valued compactly supported function such that $0\le\chi_0\le 1$ and $\chi_0=1$
on supp$\chi$. It suffices to show that the operator ${\rm Op}_h(\chi_0+b)$ is invertible. Indeed, this would imply (3.8) and (3.9) 
with $E=\left({\rm Op}_h(\chi_0+b)\right)^{-1}$ and $L_l=E{\rm Op}_h(\chi_0)$, $L_r={\rm Op}_h(\chi_0)E$.

An easy computation shows that
$$b_0=\frac{\widetilde c(x')(c_0(x')r_0(x',\xi')-z)}{c_1\rho_1+c_2\rho_2}(1-\chi(x',\xi')),$$
where $\widetilde c$ and $c_0$ are the restrictions on $\Gamma$ of the functions
$$c_1n_1-c_2n_2\quad\mbox{and}\quad\frac{c_1^2- c_2^2}{c_1n_1-c_2n_2},$$
respectively. Let us see that
$$C_1\langle\xi'\rangle^k\le|\chi_0+b_0|\le C_2\langle\xi'\rangle^k\eqno{(3.10)}$$
with some constants $C_1,C_2>0$, where $k=-1$ if $c_0(x') \equiv 0$ and
$k=1$ if $c_0(x')\neq 0$, $\forall x'\in \Gamma$. Since 
$$b_0=\frac{\widetilde c(c_0r_0-z)}{i(c_1+c_2)\sqrt{r_0}}(1-\chi)\left(1+{\mathcal O}\left(\langle\xi'\rangle^{-1}\right)\right),$$
 we have with some positive constants $\widetilde C$, $\widetilde C_1$, $\widetilde C_2$,
$$2|\chi_0+b_0|\ge |\chi_0+{\rm Re}\,b_0|+|{\rm Im}\,b_0|\ge \chi_0-|{\rm Re}\,b_0|+|{\rm Im}\,b_0|$$
$$\ge \chi_0+\frac{\widetilde C_1}{\langle\xi'\rangle}|c_0r_0-z|(1-\chi)\left(1-{\mathcal O}\left(\langle\xi'\rangle^{-1}\right)\right)$$ $$
\ge \chi_0+\widetilde C_2\langle\xi'\rangle^k(1-\chi)\ge\widetilde C\langle\xi'\rangle^k $$
which yields the lower bound in (3.10). The upper bound is obvious.

It follows from (3.10) that $(\chi_0+b_0)^{-1}\in S^{-k}$. Moreover, when $c_1\equiv c_2$, $\partial_\nu c_1\equiv \partial_\nu c_2$ on $\Gamma$, 
by (2.8) one concludes that
$b_1\in S^{-2}$. Hence the operator ${\rm Op}_h(\chi_0+b)$ is invertible with an inverse which is an $h-\Psi$DO with a symbol
belonging to the class $S^{-k}$. In particular, we have $\left({\rm Op}_h(\chi_0+b)\right)^{-1}:H_h^s\to H_h^{s+k}= {\cal O}(1)$, 
${\rm Op}_h(\chi_0+b):H_h^s\to H_h^{s-k}= {\cal O}(1)$, $\forall s\in \R$.

\eproof

Set ${\mathcal V}_j(h):=\{\nu_k\in{\rm spec}\,G_D^{(j)}:h^2\nu_k\in Z\}$, $j=1,2$. 
Define the operator ${\cal K}$ as follows:
$${\mathcal K}(z,h)=E(z,h)\left(c_1F_1(z,h)-c_2F_2(z,h)\right)+L_l(z,h)\quad{\rm if}\quad k=-1,$$
$${\mathcal K}(z,h)=\left(c_1F_1(z,h)-c_2F_2(z,h)\right)E(z,h)+L_r(z,h)\quad{\rm if}\quad k=1.$$
We obtain easily that
$$E (c_1 {\mathcal N}_1 - c_2 {\mathcal N}_2) = I + {\mathcal K}\: \quad{\rm if}\: \quad k = -1,$$
$$(c_1 {\mathcal N}_1 - c_2 {\mathcal N}_2) E = I + {\mathcal K}\:\quad {\rm if} \: \quad k = 1.$$
Clearly, the operator ${\mathcal K}$ 
is trace class and meromorphic in $z\in Z$ with poles  $\{w_k\}$, $w_k/h^2\in{\cal V}_1(h)\cup{\cal V}_2(h)$, and residue of finite rank, 
so we can define the meromorphic function
$$g_h(z):={\rm det}\left(I+{\mathcal K}(z,h)\right).$$

\begin{lemma} For all $z\in Z$ such that 
$$\delta^\sharp(z,h):=\min\{1,{\rm dist}\{z,{\rm spec}\,h^2G_D^{(1)}\cup {\rm spec}\,h^2G_D^{(2)}\}\}>0$$
we have the bound
$$\log|g_h(z)|\le C_\varepsilon h^{1-d}\delta^\sharp(z,h)^{-\varepsilon},\quad\forall\,0<\varepsilon\ll 1.\eqno{(3.11)}$$
\end{lemma}

{\it Proof.} It follows from Lemma 2.2 and the properties of the characteristic values that $\mu_j({\cal K})$ satisfy the bound
(2.10) with a new constant $C>0$ and $\delta$ replaced by $\delta^\sharp$. In fact, for $k = -1$ we have
$$\mu_{j_1 + j_2-1}({\mathcal K}) \leq \mu_{j_1}\Bigl((E(z,h)\left(c_1F_1(z,h)-c_2F_2(z,h)\right)\Bigr) + \mu_{j_2}(L_l(z, h)).$$
Since the operator $E(x, h)$ is bounded, for the first term on the right hand side we apply Lemma 2.2. On the other hand, 
$\mu_{j_2}(E Op_h(\chi_0)) \leq C \mu_{j_2} (Op_h(\chi_0))$ and for $\mu_{j_2}(Op_h(\chi_0))$ we obtain easily (2.10) with $\delta(z, h) = 1$ since $\chi_0$ has compact support. Next, if $j = j_1 + j_2 - 1$, then we have $j_1 \geq (j+1)/2$ or $j_2 \geq (j+1)/2.$ The case $k = 1$ is similar.\\

Therefore, we have
$$\log|g_h(z)|\le \sum_{j=1}^\infty\log\left(1+\mu_j({\cal K})\right)\le
\sum_{j=1}^\infty\log\left(1+C\delta^\sharp(z,h)^{-1}h^{-2m}j^{-2m/(d-1)}\right)$$
$$\le \int_{0}^\infty\log\left(1+C\delta^\sharp(z,h)^{-1}h^{-2m}t^{-2m/(d-1)}\right)dt$$
 $$=C_m h^{-d+1}\left(\delta^\sharp(z,h)\right)^{-\frac{d-1}{2m}}\int_{0}^\infty\log\left(1+t^{-2m/(d-1)}\right)dt$$
 $$\le {\widetilde C_m}h^{-d+1}\delta^\sharp(z,h)^{-\frac{d-1}{2m}}.\eqno{(3.12)}$$
 Now, given any $0<\varepsilon\ll 1,
$ we can take $m\sim \frac{d-1}{2\varepsilon}$ and $N\ge 4m$, and (3.11) follows 
 from (3.12).
 \eproof
 
 The next lemma is an almost direct consequence of the results of \cite{kn:V}.
 
 \begin{lemma} Let $\kappa$ be as in Theorem $1.1$. Then, given any $0<\epsilon\ll 1,$ the operator $I+{\cal K}(z,h)$ 
 is invertible on $L^2(\Gamma)$
 for $z\in Z$, $|{\rm Im}\,z|\ge h^{\kappa-\epsilon},$ and its inverse satisfies in this region the bound
 $$\left\|\left(I+{\cal K}(z,h)\right)^{-1}\right\|_{L^2\to L^2}\le Ch^{-\ell}\eqno{(3.13)}$$
 with some constants $C,\ell>0$. For these values of $z$ we also have
 $$\log\frac{1}{|g_h(z)|}\le C_\varepsilon h^{1-d-\varepsilon},\quad\forall\,0<\varepsilon\ll 1.\eqno{(3.14)}$$
 Moreover, the function $\log g_h(z)$ is holomorphic in $z\in Z$, $|{\rm Im}\,z|\ge h^{\kappa-\epsilon}$ and satisfies the bound
 $$\left|\frac{d}{dz}\log g_h(z)\right|\le \frac{C_\epsilon h^{1-d-2\epsilon}}{|{\rm Im}\,z|}\eqno{(3.15)}$$
 in $W:=\{z\in \C:\frac{2}{3}\le|{\rm Re}\,z|\le \frac{5}{2},\, 2h^{\kappa-\epsilon}\le |{\rm Im}\,z|\le \frac{1}{2}\}$.
 \end{lemma} 

{\it Proof.} It follows from the analysis in Section 5 of \cite{kn:V} that, under the assumptions of Theorem 1.1, 
the operator $c_1{\cal N}_1(z,h)-c_2{\cal N}_2(z,h)$ is invertible
for $z\in Z$, $|{\rm Im}\,z|\ge h^{\kappa-\epsilon}$ and
$$\left\|\left(c_1{\cal N}_1(z,h)-c_2{\cal N}_2(z,h)\right)^{-1}\right\|_{H_h^1\to L^2}\le Ch^{-\ell}\quad{\rm if}\quad k=-1,$$
$$\left\|\left(c_1{\cal N}_1(z,h)-c_2{\cal N}_2(z,h)\right)^{-1}\right\|_{L^2\to H_h^1}\le Ch^{-\ell}\quad{\rm if}\quad k=1.$$
Now (3.13) follows from these bounds and Lemma 3.3 because
$$\left(I+{\cal K}\right)^{-1}=\left(c_1{\cal N}_1-c_2{\cal N}_2\right)^{-1}E^{-1}\quad{\rm if}\quad k=-1,$$
$$\left(I+{\cal K}\right)^{-1}=E^{-1}\left(c_1{\cal N}_1-c_2{\cal N}_2\right)^{-1}\quad{\rm if}\quad k=1.$$
The bound (3.14) can be obtained in precisely the same way as (3.11) by using (3.13) and the formula
$$\frac{1}{g_h(z)}={\rm det}\Bigl(I-(I+{\cal K}(z,h))^{-1}{\cal K}(z,h)\Bigr).$$
Note that the norm $\|(I + {\cal K}(z, h))^{-1}\|$ will add a factor $h^{-\frac{l(d-1)}{2m} }$ which for 
sufficiently large $m$ yields a factor ${\cal O}(h^{-\epsilon})$.

Clearly, it follows from the Fredholm theorem that, under the assumptions of Theorem 1.1, the operator-valued
function $\left(I+{\cal K}(z,h)\right)^{-1}:L^2(\Gamma)\to L^2(\Gamma)$ is meromorphic in $Z$ with finite rank residue and holomorphic
with respect to $z\in Z$ for $|{\rm Im}\,z|\ge h^{\kappa-\epsilon}$. Therefore the functions $g_h(z)$ and $\frac{1}{g_h(z)}$ are holomorphic
in $z\in Z$, $|{\rm Im}\,z|\ge h^{\kappa-\epsilon}$, and hence so is $\log g_h(z)$. Fix an arbitrary $w\in W$. Then the function
$f(z)=\log\frac{g_h(z)}{g_h(w)}$ is holomorphic in $z\in Z$, $|{\rm Im}\,z|\ge h^{\kappa-\epsilon}$ and $f(w)=0$.
It follows from the bounds (3.11) and (3.14) that ${\rm Re}f(z)\le {\cal O}_{\epsilon}(h^{1-d-2\epsilon})$ 
for $z\in Z$, $|{\rm Im}\,z|\ge h^{\kappa-\epsilon}$
 In particular, the later estimate holds on the circle $C_w = \{z \in \C: |z-w|=\frac{|{\rm Im}\,w|}{2}\}$ since for every  
$z \in C_w$ we have $|\Im z| \geq \frac{ |\Im w|}{2}.$
Applying the Caratheodory theorem (e.g. see 5.5 in \cite{kn:T}), we get 
$$|f'(z)|= {\cal O}_{\epsilon}(h^{1-d-2\epsilon})|{\rm Im}\,w|^{-1}\:\quad {\rm for}\: \quad |z-w|\le \frac{|{\rm Im}\,w|}{3}.$$ 
 This implies (3.15) because $f'(z)=\frac{d}{dz}\log g_h(z)$. 
\eproof

Let $\gamma_0 \subset Z$ be a simple closed positively oriented curve which avoids the eigenvalues of 
$h^2G_D^{(j)}$, $j=1,2$, as well as the poles of $T(z/h^2)^{-1}$. Denote by $M_{\gamma_0}(h)$ the number of the poles, $\{\lambda_k\}$,
of $R(\lambda)$ such that $h^2\lambda_k$ are in the interior of the domain $\omega_0$ with boundary $\gamma_0$. Similarly, we denote by 
$M_{\gamma_0}^{(j)}(h)$ the number of the eigenvalues, $\{\nu_k\}$, of $G_D^{(j)}$ such that $h^2\nu_k$ are in $\omega_0$. Corollary 3.2
implies the following

\begin{lemma} We have the identity
$$M_{\gamma_0}(h)=M_{\gamma_0}^{(1)}(h)+M_{\gamma_0}^{(2)}(h)+\frac{1}{2\pi i}\int_{\gamma_0}
\frac{d}{dz}\log g_h(z)dz.\eqno{(3.16)}$$
\end{lemma}

{\it Proof.} We apply (3.6) and use the identities
$$h T(z/h^2) = E^{-1}(z, h)(I + {\mathcal K}(z, h)), \: (h T(z/h^2))^{-1} = ( I + {\cal K}(z, h))^{-1} E(z, h)$$
combined with the analyticity of $E(z, h)$ in $z$ and the following well-known formula
$${\rm tr}\,(I+{\mathcal K}(z,h))^{-1}\frac{d{\mathcal K}(z,h)}{dz}= \frac{d}{dz}\log{\rm det}(I+{\mathcal K}(z,h)).$$
The above formula for $\log{\rm det}(I+{\mathcal K}(z,h))$ is classical for finite rank  perturbations of the identity. For trace class ones this formula follows by an approximation with finite rank operators (see for example, Section 5, \cite{kn:Sj1}).
\eproof

It follows from (3.16) that $z_0\in Z\setminus {\rm spec}(h^2G_D^{(1)})\cup {\rm spec}(h^2G_D^{(2)})$ is a zero of $g_h(z)$
if and only if $z_0$ is a pole of $R(z/h^2)$ (and hence $z_0/h^2$ is an interior transmission eigenvalue) and the multiplicities coincide. 
Similarly, one can see that if $\widetilde z_0$ is a pole of $g_h(z)$ with multiplicity $\widetilde m_0$, then  
$\widetilde z_0\in {\rm spec}(h^2G_D^{(1)})\cup {\rm spec}(h^2G_D^{(2)})$ and $\widetilde m_0\le\widetilde m_1+\widetilde m_2$, 
where $\widetilde m_j$ is the multiplicity of $\widetilde z_0/h^2$
as an eigenvalue of $G_D^{(j)}$. In what follows we will use the formula (3.16) to prove the following

\begin{prop}For every $0<\epsilon\ll 1$ and $A>0$, independent of $h$, we have the asymptotics
$$I(h): = \sharp\left\{z_k,\,z_k/h^2\,\, \mbox{is (ITE)}:\,\, 1-Ah^{\kappa-\epsilon}\le |{\rm Re}\,z_k|\le 2+Ah^{\kappa-\epsilon},\,
|{\rm Im}\,z_k|\le h^{\kappa-\epsilon}\right\}$$ $$=(2^{d/2}-1)(\tau_1+\tau_2)h^{-d}+
{\cal O}_{\epsilon,A}(h^{-d+\kappa-3\epsilon}),\: 0 < h \leq h_0(\epsilon, A).\eqno{(3.17)}$$
\end{prop}
 
 {\it Proof.} We will consider only the case ${\rm Re}\,z_k>0$, since the case ${\rm Re}\,z_k<0$ is similar (and even simpler since the function
 $g_h(z)$ does not have poles in ${\rm Re}\,z<0$). Consider the points 
$w_1^\pm=1-Ah^{\kappa-\epsilon}\pm\frac{i}{3}$,
 $w_2^\pm=2+Ah^{\kappa-\epsilon}\pm\frac{i}{3}$, $\widetilde w_1^\pm=1-Ah^{\kappa-\epsilon}\pm i3h^{\kappa-\epsilon}$,
 $\widetilde w_2^\pm=2+Ah^{\kappa-\epsilon}\pm i3h^{\kappa-\epsilon}$ and set
 $$\Theta_1=\left\{z\in \C:1-2(A+1)h^{\kappa-\epsilon}\le {\rm Re}\,z\le 1+h^{\kappa-\epsilon},\,|{\rm Im}\,z|\le 4h^{\kappa-\epsilon}\right\},$$ 
 $$\Theta_2=\left\{z\in \C:2-h^{\kappa-\epsilon}\le {\rm Re}\,z\le 2+2(A+1)h^{\kappa-\epsilon},\,|{\rm Im}\,z|\le 4h^{\kappa-\epsilon}\right\}.$$
 The following lemma will be proved later on.
 
 \begin{lemma} There exist positively oriented piecewise smooth curves 
 $\widetilde\gamma_1\subset\Theta_1$ and  
 $\widetilde\gamma_2\subset\Theta_2$, where
 $\widetilde\gamma_1$ connects the point $\widetilde w_1^-$ with $\widetilde w_1^+$, while $\widetilde\gamma_2$ connects the point $\widetilde 
 w_2^+$ with $\widetilde w_2^-$, such that
 $$\left|{\rm Im}\,\int_{\widetilde\gamma_j}
\frac{d}{dz}\log g_h(z)dz\right|\le C_\epsilon h^{-d+\kappa-2\epsilon}, \:\quad j = 1,2.\eqno{(3.18)}$$
 \end{lemma}
 
 Now we apply Lemma 3.6 with a contour $\gamma_0 =\gamma_1\cup\gamma_3\cup\gamma_2\cup\gamma_4$, where $\gamma_3\subset W$ 
 is the segment $[w_1^+,w_2^+]$ on the line passing through the points $w_1^+$ and $w_2^+$, and 
 $\gamma_4\subset W$ is the segment $[w_2^-,w_1^-]$ on the line passing through the points $w_2^-$ and $w_1^-$.
  Next, $\gamma_1=[w_1^-,\widetilde w_1^-]\cup\widetilde\gamma_1\cup[\widetilde w_1^+,w_1^+]$, 
 $\gamma_2=[w_2^+,\widetilde w_2^+]\cup\widetilde\gamma_2\cup[\widetilde w_2^-,w_2^-]$ (see Figure 1).

\begin{figure}[tbp] % float placement: (h)ere, page (t)op, page (b)ottom, other (p)age
  \centering
  % file name: C:/Users/petkov/Documents/My PCTeX Files/vodev/contour.eps
  \includegraphics[bb=0 0 512 285,width=5.67in,height=3.15in,keepaspectratio]{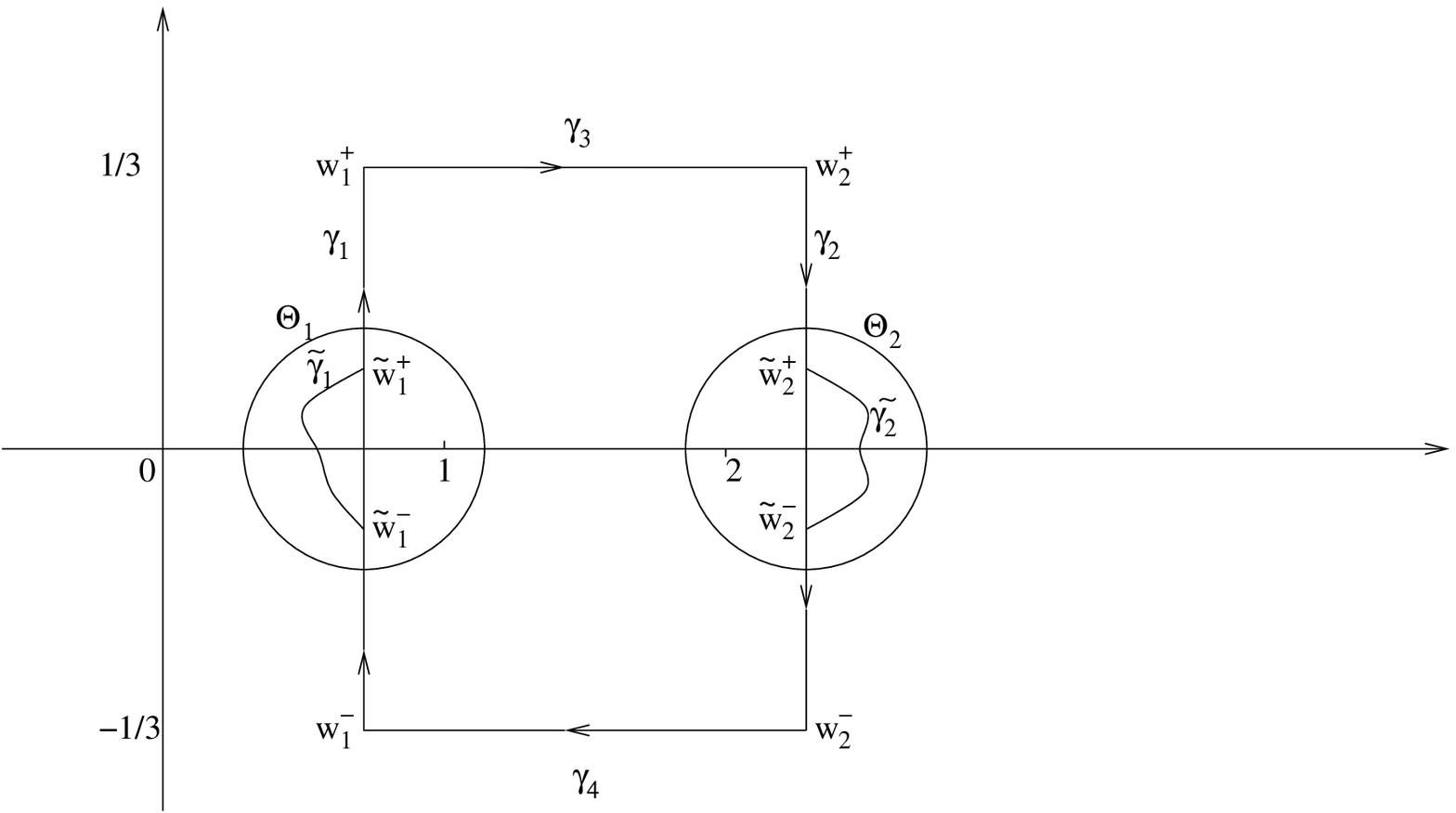}
  \caption{Contour $\gamma_0$}
  \label{fig:contour}
\end{figure}

 Since $\gamma_j\subset W$, $|\gamma_j|= {\cal O}(1)$, $j=3,4$, by (3.15) we have
 $$\left|\int_{\gamma_j}
\frac{d}{dz}\log g_h(z)dz\right|\le \int_{\gamma_j}\left|
\frac{d}{dz}\log g_h(z)\right||dz|$$ $$\le
C_\epsilon h^{-d+1-2\epsilon}\int_{\gamma_j}|dz|\le
C_\epsilon h^{-d+1-2\epsilon},\quad j = 3, 4.\eqno{(3.19)}$$
Applying (3.15) once more, we have
$$\left|\int_{[w_j^\pm,\widetilde w_j^\pm]}
\frac{d}{dz}\log g_h(z)dz\right|\le C_\epsilon h^{-d+1-2\epsilon}\int_{3h^{\kappa-\epsilon}}^{1/2}\frac{d\sigma}{\sigma}\le 
C_\epsilon h^{-d+1-3\epsilon},\:\quad j =1,2.\eqno{(3.20)}$$
On the other hand, since the counting function of the eigenvalues of $G_D^{(j)}$ satisfies the Weyl law, we deduce
$$M_{\gamma_0}^{(j)}(h)\le\sharp\left\{\nu_k\in{\rm spec}\,G_D^{(j)}:1-2(A+1)h^{\kappa-\epsilon}\le h^2\nu_k\le 2+
2(A+1)h^{\kappa-\epsilon}\right\}$$
 $$=\tau_j\left(\frac{2}{h^2}+\frac{2(A+1)h^{\kappa-\epsilon}}{h^2}\right)^{d/2}-\tau_j\left(\frac{1}{h^2}-
 \frac{2(A+1)h^{\kappa-\epsilon}}{h^2}\right)^{d/2}+{\cal O}_{\epsilon}(h^{-d+1})$$
 $$=(2^{d/2}-1)\tau_jh^{-d}+{\cal O}_{\epsilon,A}(h^{-d+\kappa-\epsilon})$$
 and similarly
 $$M_{\gamma_0}^{(j)}(h)\ge\sharp\left\{\nu_k\in{\rm spec}\,G_D^{(j)}:1+h^{\kappa-\epsilon}\le h^2\nu_k\le 2-h^{\kappa-\epsilon}\right\}$$
 $$=(2^{d/2}-1)\tau_jh^{-d}-{\cal O}_{\epsilon}(h^{-d+\kappa-\epsilon}).$$
 Consequently,
 $$M_{\gamma_0}^{(j)}(h)=(2^{d/2}-1)\tau_jh^{-d}+{\cal O}_{\epsilon,A}(h^{-d+\kappa-\epsilon}).\eqno{(3.21)}$$
Taking together (3.16), (3.18), (3.19), (3.20) and (3.21), we obtain
 $$M_{\gamma_0}(h)=(2^{d/2}-1)(\tau_1+\tau_2)h^{-d}+{\cal O}_{\epsilon,A}(h^{-d+\kappa-3\epsilon}).\eqno{(3.22)}$$
 Thus, to establish (3.17), it remains to show that the counting function $I(h)$ satisfies
 $$\left|I(h)-M_{\gamma_0}(h)\right|\le C_{\epsilon,A}h^{-d+\kappa-3\epsilon}.\eqno{(3.23)}$$
 Given a parameter $\theta > 0$, independent of $h$, introduce $B_j^\pm(\theta)=\{z\in \C:|z-\widetilde w_j^\pm|\le \theta h^{\kappa-\epsilon}\}$.
 Clearly, there exists $\theta_0>0$ such that $\Theta_j\subset B^+_j(\theta)\cup B^-_j(\theta)$, $\forall\theta\ge\theta_0,\: j = 1,2$. 
 Let $\left\{z_k^{\pm,j}\right\}$ be the zeros (repeated with their multiplicities) of $g_h(z)$ in $B_j^\pm(2\theta_0)$ and let 
 $\left\{y_k^{\pm,j}\right\}$ be the poles (repeated with their multiplicities) of $g_h(z)$ in $B_j^\pm(4\theta_0)$.
 Therefore the function
 $$f_h^{\pm,j}(z)=g_h(z)\prod_{k}\left(z-y_k^{\pm,j}\right)$$
 is holomorphic in the interior of $B_j^\pm(4\theta_0)$. Obviously, $\left\{y_k^{\pm,j}\right\}$ are among the eigenvalues of the operators $G_D^{(1)}$ and $G_D^{(2)}$
 in an interval of the form 
$$[1-{\cal O}(h^{\kappa-\epsilon}), 1+{\cal O}(h^{\kappa-\epsilon})]\cup [2-{\cal O}(h^{\kappa-\epsilon}), 
 2+{\cal O}(h^{\kappa-\epsilon})].$$
 Hence, by the Weyl law for the counting function of these eigenvalues, as in the proof of (3.21), we get
 $$\sharp\left\{y_k^{\pm,j}\right\}\le C_{\epsilon,A}h^{-d+\kappa-\epsilon},\: j = 1,2.\eqno{(3.24)}$$
  By (3.14) and (3.24), we have
  $$\log\left|f_h^{\pm,j}(\widetilde w_j^\pm)\right|=\log\left|g_h(\widetilde w_j^\pm)\right|+
  \sum_k\log\left|\widetilde w_j^\pm- y_k^{\pm,j}\right|$$
  $$\ge -C_\epsilon h^{-d+1-\epsilon}-\sharp\left\{y_k^{\pm,j}\right\}C\log\frac{1}{h}
  \ge -2C_{\epsilon, A} h^{-d+\kappa-2\epsilon}.\eqno{(3.25)}$$
  On the other hand, applying (3.11) and (3.24), for $z\in B_j^\pm(\theta), \theta_0< \theta < 4 \theta_0$, $\left|z-y_k^{\pm,j}\right|\ge h^M$, 
  $M\gg 1$, we obtain
  $$\log\left|f_h^{\pm,j}(z)\right|=\log\left|g_h(z)\right|+\sum_k\log\left|z-y_k^{\pm,j}\right|$$
  $$\le C_\epsilon h^{-d+1-\epsilon}+\sharp\left\{y_k^{\pm,j}\right\}M\log\frac{1}{h}
  \le 2C_{\epsilon, A} h^{-d+\kappa-2\epsilon}.\eqno{(3.26)}$$
 We claim that there exists $3\theta_0<\mu_1<4\theta_0$ such that the distance between $\left\{y_k^{\pm,j}\right\}$
 and the circle $\partial B^\pm_j(\mu_1)$ is greater than $h^M$, provided $M\gg d$. Indeed, if we suppose the contrary, this would imply
 that the length of the interval $J_j^\pm:=\R\cap\left(B_j^\pm(4\theta_0)\setminus B_j^\pm(3\theta_0)\right)$
 is upper bounded by $\sharp\left\{y_k^{\pm,j}\in J_j^{\pm}\right\} h^M={\cal O}(h^{M-d})$, which is impossible if $M$
 is taken large enough. This proves the claim. Thus, by (3.26) we have the estimate 
 $\log\left|f_h^{\pm,j}(z)\right|={\cal O}_{\epsilon}(h^{-d+\kappa-2\epsilon})$ on 
 $\partial B^\pm_j(\mu_1)$, which in turn implies $\log\left|f_h^{\pm,j}(z)\right|={\cal O}_{\epsilon}(h^{-d+\kappa-2\epsilon})$ on 
 $B^\pm_j(3\theta_0)$. Combining this with (3.25) and the Jensen theorem (see for example 3.6 in \cite{kn:T}), yields for the zeros $z_k^{\pm, j}$ in $B^{\pm}_j( 2 \theta_0)$ the following bound
 $$\sharp\left\{z_k^{\pm,j}: \: z_k^{\pm, j} \in B^{\pm}_j(2\theta_0)\right\}\le C_{\epsilon,A} h^{-d+\kappa-2\epsilon}.\eqno{(3.27)}$$
 Since the left-hand side of (3.23) is upper bounded by the number of the zeros and the poles of the function $g_h(z)$ in
 $B_1^+(\theta_0)\cup B_1^-(\theta_0)\cup B_2^+(\theta_0)\cup B_2^-(\theta_0)$, the estimate (3.23) follows from (3.24) and (3.27).
 \eproof
 
 \begin{rem} The bound $(3.27)$ of the number of the zeros $z_k^{\pm, j}$ of $g_h(z)$ in $B^{\pm}_j(2\theta_0)$ does 
 not depend on the statement of Lemma $3.8$ but only on the application of the Jensen theorem based on $(3.25)$, $(3.26)$. 
 We will use $(3.27)$ in the proof of Lemma $3.8$ below. 
\end{rem}

 {\it Proof of Lemma $3.8$.} We will consider only the case $j=1$, since the case $j=2$ is similar. Introduce the function
 $$\zeta_h(z): =g_h(z)\prod_{w\in{\cal M}_1}\left(z-w\right)^{-1}
 \prod_{w\in{\cal M}_2}\left(z-w\right),$$
 where ${\cal M}_1=\{z_k^{+,1}\}\cup\{z_k^{-,1}\}$ is the set of all zeros of $g_h(z)$ in $B_1^-(2\theta_0)\cup B_1^+(2\theta_0)$
 and ${\cal M}_2=\{y_k^{+,1}\}\cup\{y_k^{-,1}\}$ is the set of all poles of $g_h(z)$ in $B_1^-(4\theta_0)\cup B_1^+(4\theta_0)$.
 Since $\zeta_h(z)$ does not have zeros and poles in $B_1^-(2\theta_0)\cup B_1^+(2\theta_0)$,
  the function $\log \zeta_h(z)$ is holomorphic in $B_1^-(2\theta_0)\cup B_1^+(2\theta_0)$.
  We need the following
  
  \begin{lemma} The function $\zeta_h(z)$ satisfies the bound
  $$\log|\zeta_h(z)|\le C_{\epsilon}h^{-d+1-2\epsilon},\quad\forall z\in B_1^-(\theta)\cup B_1^+(\theta),\eqno{(3.28)}$$
 for every $0<\theta< 2\theta_0$ independent of $h$.
    \end{lemma}
 {\it Proof.} Set $U=\cup_{w\in{\cal M}}\{z\in \C:|z-w|\le h^M\}$, where $M\gg d$ and ${\cal M}={\cal M}_1\cup
 {\cal M}_2$. Clearly, $U=\cup_\nu U_\nu$, where every $U_\nu$ is a domain with a piecewise smooth boundary and $U_\nu\cap U_\mu=\emptyset$
 if $\nu\neq\mu$. Moreover, we have
 $$\sum_\nu {\rm measure} \:(\partial U_\nu)\le 2\pi h^M\,\sharp\{w\in{\cal M}\}\le Ch^{M-d}.$$
 Let $\theta<\theta_1<2\theta_0$ be independent of $h$. 
 Let $\{U^\pm_{\nu_i}\}$ be the set of all $U_\nu$ such that $U_\nu\cap\partial B_1^\pm(\theta_1) \neq \emptyset.$
 We now construct a closed curve, $\beta^\pm_1(\theta_1)$ as follows: we keep all arcs on 
 $\partial B_1^\pm(\theta_1)$ having no common points with $\{U^\pm_{\nu_i}\}$ and replace the arc  
$ \partial B_1^\pm(\theta_1)\cap U^\pm_{\nu_i}$ with arcs on $\partial U^\pm_{\nu_i}$ connecting the corresponding end points. 
Thus we can guarantee that $\beta_1^\pm(\theta_1)$ belongs to an ${\cal O}(h^{M-d})$ neighborhood
of ${\partial B_1^\pm}(\theta_1)$ and, moreover,  the distance between $\beta_1^\pm(\theta_1)$ and 
the set ${\cal M}$ is greater than $h^M$.
 In the same way, as in the proof of (3.25) and (3.26) above, by using (3.11), (3.14), (3.24) and (3.27), we get
 $$\log|\zeta_h(z)|\le C_{\epsilon}h^{-d+1-2\epsilon},\quad\forall z\in\ \beta_1^\pm(\theta_1).\eqno{(3.29)}$$
 Since $B_1^\pm(\theta)$ is in the interior of the domain bounded by $\beta_1^\pm(\theta_1)$, 
 the estimate (3.29) implies (3.28).
  \eproof
  
  We will now construct the curve  $\widetilde\gamma_1$.  Let $\{U_{\nu_i}\}$ be the set of all 
  $U_\nu$ such that $U_\nu\cap [\widetilde w_1^-,\widetilde w_1^+] \neq \emptyset.$ We keep all segments on
  $[\widetilde w_1^-,\widetilde w_1^+]$ having no common points with $\{U_{\nu_i}\}$ and replace the segments on
  $[\widetilde w_1^-,\widetilde w_1^+]\cap U_{\nu_i}$ with arcs on $\partial U_{\nu_i}$ connecting the corresponding end points.
  Thus we get a piecewise smooth curve $\widetilde\gamma_1$ belonging to an ${\cal O}(h^{M-d})$ neighborhood
of $[\widetilde w_1^-,\widetilde w_1^+]$ and the distance between $\widetilde\gamma_1$ and the set ${\cal M}$ is greater than $ h^M$.
  Hence $\widetilde\gamma_1\subset\Theta_1$.  Now we can write
 $$\int_{\widetilde\gamma_1}\frac{d}{dz}\log g_h(z)dz=
 \int_{[\widetilde w_1^-,\widetilde w_1^+]}\frac{d}{dz}\log\zeta_h(z)dz$$
 $$+\sum_{w\in{\cal M}_1}\int_{\widetilde\gamma_1}(z-w)^{-1}dz -\sum_{w\in{\cal M}_2}\int_{\widetilde\gamma_1}(z-w)^{-1}dz.\eqno{(3.30)}$$
 We will show that 
 $$\left|\frac{d}{dz}\log\zeta_h(z)\right|\le C_{\epsilon}h^{-d+1-\kappa-\epsilon},\quad\forall z\in\Theta_1,\eqno{(3.31)}$$
 $$\left|{\rm Im}\,\int_{\widetilde\gamma_1}(z-w)^{-1}dz\right|\le 3\pi,\quad\forall w\in{\cal M}.\eqno{(3.32)}$$
 Since the length of the interval $|\widetilde w_1^-,\widetilde w_1^+]$ is $6 h^{\kappa-\epsilon}$, the estimate (3.31) implies that
 the absolute value of the first integral on the right-hand side of (3.30) is ${\cal O}_{\epsilon}(h^{-d+1-2\epsilon})$. 
 Thus, (3.18) would follow
 from (3.30), (3.31), (3.32) and the bounds (3.24) and (3.27).
 
 To prove (3.31), we apply the Caratheodory theorem (see 5.5. \cite{kn:T}) for the derivative of the function $f_\pm(z)=\log\frac{\zeta_h(z)}{\zeta_h(\widetilde w_1^\pm)}$.
 Note that $\log|\zeta_h(\widetilde w_1^\pm)|$ can be bounded from below in the same way as in (3.25) above. Therefore, applying (3.28),
 we get for the real part of $f_{\pm}(z)$ the estimate
 $${\rm Re}\,f_\pm(z)= \log|\zeta_h(z)|-\log|\zeta_h(\widetilde w_1^\pm)|\le 
 Ch^{-d+1-2\epsilon},\quad\forall z\in \partial B_1^\pm\Bigl(\frac{3}{2}\theta_0\Bigr).$$
 Since $f_\pm(\widetilde w_1^\pm)=0$, we conclude by the Caratheodory theorem that $|f'_\pm(z)|= {\cal O}_{\epsilon}(h^{-d+1-\kappa-\epsilon})$ in the disc $B_1^\pm(\theta_0)$, which clearly implies (3.31).
 
 To establish (3.32), observe that if $w$ does not lie on the line connecting the points $\widetilde w_1^-$ and 
$ \widetilde w_1^+$ and if $\sigma_0>0$ denotes the distance from $w$ to this line, after a suitable change of variables,
we have 
$$\left|{\rm Im}\,\int_{\widetilde w_1^-}^{\widetilde w_1^+}(z-w)^{-1}dz\right|=\int_a^b\frac{\sigma_0d\sigma}{\sigma_0^2+\sigma^2}\le
\int_{-\infty}^\infty\frac{d\sigma}{1+\sigma^2} = \pi.\eqno{(3.33)}$$
Since the integral in the left-hand side of (3.32) differs from the integral in the left-hand side of (3.33) either by $0$ or $2\pi i$, 
the estimate (3.33) implies (3.32) in this case. If $w$ lies on the line connecting the points $\widetilde w_1^-$ and 
$ \widetilde w_1^+$, then the integral on the left-hand side of (3.32) is a limit of integrals of the first kind, and hence
(3.32) will be true in this case, too. This completes the proof of Lemma 3.8.
\eproof

{\it Proof of Theorem $1.1$}. Let $\kappa$ be as described in Theorem 1.1. 
Let $A_1$ and $A_2$ be arbitrary real numbers, independent of $h$, and let $A >\max\{|A_1|,|A_2|\}$ be independent of $h$.
It follows from the proof of Proposition 3.7 (see (3.27)) that
$$\sharp\left\{z_k,\,z_k/h^2\,\, \mbox{is (ITE)}:\,\, 1-Ah^{\kappa-\epsilon}\le |{\rm Re}\,z_k|\le 1+Ah^{\kappa-\epsilon},\,
|{\rm Im}\,z_k|\le {\cal O}(h^{\kappa-\epsilon})\right\}$$ $$=
{\cal O}_{\epsilon,A}(h^{-d+\kappa-2\epsilon}),$$
$$\sharp\left\{z_k,\,z_k/h^2\,\, \mbox{is (ITE)}:\,\, 2-Ah^{\kappa-\epsilon}\le |{\rm Re}\,z_k|\le 2+Ah^{\kappa-\epsilon},\,
|{\rm Im}\,z_k|\le {\cal O}(h^{\kappa-\epsilon})\right\}$$ $$=
{\cal O}_{\epsilon,A}(h^{-d+\kappa-2\epsilon}).$$
Therefore, by (3.17) we get for every $0 < \epsilon \ll 1$
$$\sharp\left\{z_k,\,z_k/h^2\,\, \mbox{is (ITE)}:\,\, 1-A_1h^{\kappa-\epsilon}\le |{\rm Re}\,z_k|\le 2+A_2h^{\kappa-\epsilon},\,
|{\rm Im}\,z_k|\le {\cal O}(h^{\kappa-\epsilon})\right\}$$ $$=(2^{d/2}-1)(\tau_1+\tau_2)h^{-d}+
{\cal O}_{\epsilon,A_1,A_2}(h^{-d+\kappa-3\epsilon}), \: 0 < h \leq h_1(A_1, A_2, \epsilon).$$
Choose $h = \frac{\sqrt{2}}{r},\: r \gg 1.$ The above asymptotics  
yields
$$\left\{\lambda\in\C:\,\lambda \: {\text is\: (ITE)},\: \frac{r^2}{2}-A_1r^{2-\kappa+\epsilon}\le|{\rm Re}\,\lambda|\le r^2+
A_2r^{2-\kappa+\epsilon},\,|{\rm Im}\,\lambda|\le r^{2-\kappa+\epsilon}\right\}$$
$$= (1- 2^{-d/2})(\tau_1+\tau_2)r^{d}+{\cal O}_{\epsilon, A_1,A_2}(r^{d -\kappa +3\epsilon}),\quad r \geq r_1(A_1,A_2,\epsilon).$$
 Recall that according to the results in \cite{kn:V}, there are no (ITE) in the region 
$$\left\{\lambda\in\C:\, \frac{r^2}{2}\le|\lambda|\le r^2,\,\,|{\rm Im}\,\lambda|\ge r^{2-\kappa+\epsilon}\right\}$$
 for every 
$0<\epsilon\ll 1$, provided $r\ge r_0(\epsilon)\gg 1$. On the other hand, it is clear that the region
 $$\left\{\lambda\in\C:\, \frac{r^2}{2}\le| \lambda|\le r^2,\,\,|{\rm Im}\,\lambda|\le r^{2-\kappa+\epsilon}\right\}$$
 is contained in the region 
 $$\left\{\lambda\in\C:\: \frac{r^2}{2}-r^{2-\kappa+\epsilon}\le|{\rm Re}\,\lambda|\le r^2,\,\,|{\rm Im}\,\lambda|\le r^{2-\kappa+\epsilon}\right\}$$
and contains the region
$$\left\{\lambda\in\C:\: \frac{r^2}{2}\le|{\rm Re}\,\lambda|\le r^2-
r^{2-\kappa+\epsilon},\,\,|{\rm Im}\,\lambda|\le r^{2-\kappa+\epsilon}\right\}.$$
Thus we get the asymptotics
$$N(r)-N(r/\sqrt{2})=(1-2^{-d/2})(\tau_1+\tau_2)r^{d}+{\cal O}_{\varepsilon}(r^{d-\kappa+\varepsilon}), \: \quad r \geq r_0(\epsilon),\eqno{(3.34)}$$
for every $0<\varepsilon\ll 1$. Here we replace $3\epsilon$ by $\epsilon$, which is not important since our argument works for every
$0 < \epsilon \ll 1.$ The asymptotics (3.34) yields
$$N(r/2^{k/2})-N(r/2^{(k+1)/2})=(2^{-kd/2}-2^{-(k+1)d/2})(\tau_1+\tau_2)r^{d}+2^{-kd/2}{\cal O}_{\varepsilon}(r^{d-\kappa+\varepsilon})\eqno{(3.35})$$
for every integer $k\ge 0$ such that $r 2^{-k/2}\geq r_0(\epsilon)$. Let $k_0(r) \in \N$ be the smallest integer 
such that $r 2^{-k_0(r)/2} < r_0(\epsilon)$. It is clear that we have
$$N(r/2^{(k_0(r)+1)/2}) \leq N(r_0(\epsilon))=R_0(\epsilon) \eqno{(3.36)}$$
 with a constant $R_0(\epsilon) > 0$ independent of $r.$ Moreover,
$$\Bigl(2^{-(k_0(r) + 1)/2} r\Bigr)^d \leq (r_0(\epsilon))^d = R_1(\epsilon)$$
with $R_1(\epsilon) > 0$ independent of $r$.  Summing up the asymptotics (3.35) and using (3.36), we get (1.7). Thus the proof of Theorem 1.1 is complete.
\eproof

\end{document}